\documentclass[12pt]{article}
\usepackage[T1]{fontenc}
\usepackage[latin9]{inputenc}
\usepackage{color,soul}
\usepackage{framed}
\usepackage{amsmath}
\usepackage{amsthm}
\usepackage{amssymb}

\makeatletter

\providecommand{\tabularnewline}{\\}

\numberwithin{equation}{section}
\numberwithin{figure}{section}
\theoremstyle{plain}

  \theoremstyle{remark}

\usepackage{fullpage}\usepackage{makeidx}\usepackage{amsfonts}\usepackage{latexsym}

\def\C{{\mathbb{C}}}
\def\Z{{\mathbb{Z}}}
\def\R{{\mathbb{R}}}

\def\G{{\mathcal{G}}}
\DeclareMathOperator{\Hom}{Hom}
\theoremstyle{definition}
\newtheorem{lemma}{Lemma}[section]
\newtheorem{theorem}[lemma]{Theorem}\newtheorem{corollary}[lemma]{Corollary}\newtheorem{proposition}[lemma]{Proposition}\newtheorem{definition}[lemma]{Definition}\newtheorem{remark}[lemma]{Remark}\usepackage{times}

\usepackage{babel}

\usepackage[blocks]{authblk}
\title{Vertex operator algebras generated by two Ising vectors}

\author{Xiangyu Jiao\footnote{The author is supported by the NSFC No.~11971167 and is partially supported
by the STCSM No.~18dz2271000.}}
\affil{School of Mathematical Science, East China Normal University, Shanghai 200241, China}
\author{Wen Zheng}
\affil{School of Mathematics and Statistics, Qingdao University, Qingdao 266071, China}

\usepackage{ulem}
\usepackage{cancel}

\usepackage{babel}

\usepackage{babel}

\usepackage{babel}

\usepackage{babel}
\providecommand{\claimname}{Claim}
\providecommand{\theoremname}{Theorem}

\usepackage{babel}
\providecommand{\claimname}{Claim}
\providecommand{\theoremname}{Theorem}

\numberwithin{equation}{section}

\usepackage{babel}
\providecommand{\claimname}{Claim}
\providecommand{\theoremname}{Theorem}

\makeatother
\usepackage{babel}
  \providecommand{\claimname}{Claim}
\providecommand{\theoremname}{Theorem}
\def\ha{\frac{1}{2}}
\theoremstyle{remark}

\begin{document}
\definecolor{shadecolor}{rgb}{0,1,0}
\maketitle
\begin{abstract}
 In this paper we mainly study the vertex operator algebra $\C \mathrm{VA}(e, f)$ generated by two Ising vectors $e$ and $f$ with $\langle e, f\rangle=\frac{5}{2^{10}}$. We prove $\C \mathrm{VA}(e, f)$ is isomorphic to the $6A$-algebra $\mathcal{U}_{6A}$ constructed in \cite{LYY}. We also discuss the cases $\langle e, f\rangle=\frac{1}{4}, \frac{1}{2^{5}}, 0 $ and the case $\langle e, f\rangle=\frac{1}{2^{8}}$ with $(\tau_{e}\tau_{f})^4=1$.
\end{abstract}
\section{Introduction}
\sethlcolor{green}

In the study of the moonshine vertex operator algebra (VOA) $V^{\natural}$ constructed in \cite{FLM}, Dong, et al. in \cite{DMZ} showed that
$V^{\natural}$ contains $48$ Virasoro vectors, each Virasoro vector generates a Virasoro vertex operator algebra isomorphic to $L(\ha,0)$ in $V^{\natural}$ and $L(\ha,0)^{\otimes 48}$ is a conformal subalgebra of $V^{\natural}.$ Such a Virasoro vector is called an Ising vector. Later, Miyamoto in \cite{M1} constructed a $\tau$-involution $\tau_{e}$ for each Ising vector $e$ and showed that each axis of the monstrous Griess algebra in
\cite{C} is essentially a half of an Ising vector $e$ and $\tau_{e}$ is a $2A$-involution of the Monster simple group $\mathbb{M}$ constructed by Griess
\cite{G}. It was proved in \cite{C} that the conjugacy class of the product of two $2A$-involutions $\tau\tau'$ is one of the nine classes $1A, 2A, 3A, 4A, 5A, 6A, 4B, 2B$ and $3C$ in  $\mathbb{M}$ and the inner product of the axis $e_{\tau}$, $e_{\tau'}$ is uniquely determined by the conjugacy
class. The above result is listed in terms of Ising vector and $\tau$-involution as follows:
\begin{center}
\begin{tabular}{|c|c|c|c|c|c|c|c|c|c|}
\hline
$\left\langle \tau_{e}\tau_{f}\right\rangle ^{\mathbb{M}}$  & $1A$  & $2A$  & $3A$  & $4A$  & $5A$  & $6A$  & $4B$  & $2B$  & $3C$\tabularnewline
\hline
$\left\langle e,f\right\rangle $  & $1/4$  & $1/2^{5}$  & $13/2^{10}$  & $1/2^{7}$  & $3/2^{9}$  & $5/2^{10}$  & $1/2^{8}$  & $0$  & $1/2^{8}$\tabularnewline
\hline
\end{tabular}
\par\end{center}

A VOA $V=\oplus_{n=0}^{\infty}V_{n}$ is of \emph{moonshine type} if $\dim V_{0}=1$ and $V_{1}=0$. The inner product of any two Ising vectors of a moonshine
type VOA was determined in \cite{M2} \cite{S} and these inner products are just given in the table. Lam, et al. in \cite{LYY}, \cite{LYY1} constructed nine coset subalgebras $\mathcal{U}_{nX}$ of the lattice vertex operator algebra $V_{\sqrt{2}E_{8}}$ corresponding to the type $nX$ of $\left\langle \tau_{e}\tau_{f}\right\rangle ^{\mathbb{M}}$. Moreover, each $\mathcal{U}_{nX}$ is generated by two Ising vectors whose inner product is exactly those given in the table.  It is natural to ask whether
the subVOA generated by any two Ising vectors in a moonshine type VOA is isomorphic to one of $\mathcal{U}_{nX}$. For the cases $\langle e, f\rangle=\frac{13}{2^{10}}$, $\langle e, f\rangle=\frac{3}{2^{9}}$ and $\langle e, f\rangle=\frac{1}{2^{8}}$ with $|\tau_{e}\tau_{f}|=3$,  the answer is yes \cite{M2} \cite{SY}\cite{Zh}. In this paper, we give a positive answer for all other cases except the case $\langle e, f \rangle=\frac{1}{2^{7}}$.

We now give a brief review of our main ideas. Let $V=\oplus_{n=0}^{\infty}V_{n}$ be a moonshine type VOA and $e, f$ be two Ising vectors in $V$. Let
$\C\mathrm{VA}(e, f)$ be the subVOA generated by $e$ and $f$. First we study the Griess algebra $\C\mathrm{VA}(e, f)_{2}$. We show that $ \C\mathrm{VA}(e, f)_{2}$ always contains a set of mutually orthogonal conformal vectors such that their sum is the Virasoro element of $\C\mathrm{VA}(e, f)$
and the central charges of these conformal vectors are all coming from the unitary series
 $$
c=c_{m}=1-\frac{6}{(m+2)(m+3)}$$
where $m\in \{1,2,3,...\}$. Such a conformal vector generates a simple Virasoro VOA isomorphic to $L(c_{m}, 0)$ inside
$\C\mathrm{VA}(e, f)$. Then we use representation theory of $L(c_{m}, 0)$ to determine the structure of  $\C\mathrm{VA}(e, f)$. Though we used similar ideas
as in \cite{Zh}, the arguments here for the case $\langle e, f\rangle=\frac{5}{2^{10}}$ are more complicated. Unfortunately, the idea we used in this paper does not work for the case $\langle e, f\rangle=\frac{1}{2^{7}}$ since in this case we cannot decompose the Virasoro element into a sum of mutually orthogonal conformal vectors  coming from the minimal series. So we will not discuss it in this paper.

The paper is organized as follows. In Section 2, we review some basic
notions and some results needed later. Section 3 is the main part. First we study the case $\langle e, f \rangle=\frac{5}{2^{10}}$ and prove that in this case the VOA $\C\mathrm{VA}(e, f)$ is isomorphic to the $6A$-algrebra $\mathcal{U}_{6A}$. Then we study other cases.

\section{Preliminary}

This section is mostly expository. We recall basic definitions of various notions on vertex operator algebras. We also give a brief review of some results about  quantum dimensions and quantum Galois theory from \cite{DM2},\cite{DLM3},\cite{DJX}, the unitary series \cite{DMZ}, \cite{W} of the Virasoro algebras. In addition, we discuss the classification of the Griess algebras generated by two Ising vectors \cite{S}.

\subsection{Basics}
Let $(V,Y,\mathbf{1},\omega)$ be a vertex operator algebra and let $Y(v,z)=\sum_{n\in \Z}v_n z^{-n-1}$ denote the vertex operator of $V$ for $v\in V$, see \cite{B}.

\begin{definition}
A vertex operator algebra $V$ is said to be \emph{CFT type} if $V=\oplus_{n\in\Z_{+}}V_{n}$ and $V_{0}=\C\mathbf{1}$.
\end{definition}

\begin{definition}
An \emph{automorphism $g$} of a vertex operator algebra $V$ is a linear isomorphism of $V$ satisfying $g(\omega)=\omega$ and $gY(v, z)g^{-1}=Y(gv, z)$ for any $v\in V$. We denote by $\mathrm{Aut}(V)$ the group of all automorphisms of $V$.
\end{definition}

For a subgroup $G\le\mbox{Aut}\left(V\right)$ the fixed point set
$V^{G}=\left\{ v\in V|g\left(v\right)=v,\forall g\in G\right\} $
has a vertex operator algebra structure.

Let $g$ be an automorphism of a vertex operator algebra $V$ of order
$T$. Denote the decomposition of $V$ into eigenspaces of $g$ as:

\[
V=\oplus_{r\in\mathbb{Z}/T\text{\ensuremath{\mathbb{Z}}}}V^{r}
\]
where $V^{r}=\left\{ v\in V|gv=e^{2\pi ir/T}v\right\} $.

\begin{definition}
A \emph{weak $g$-twisted $V$-module} $M$ is
a vector space with a linear map
\begin{align*}
Y_{M}: & V\to\left(\text{End}M\right)\{z\}\\
 & v\mapsto Y_{M}\left(v,z\right)=\sum_{n\in\mathbb{Q}}v_{n}z^{-n-1}\ \left(v_{n}\in\mbox{End}M\right)
\end{align*}
which satisfies the following: for all $0\le r\le T-1$, $u\in V^{r}$,
$v\in V$, $w\in M$,
\[
Y_{M}\left(u,z\right)=\sum_{n\in-\frac{r}{T}+\mathbb{Z}}u_{n}z^{-n-1},
\]

\[
u_{l}w=0\ for\ l\gg0,
\]

\[
Y_{M}\left(\mathbf{1},z\right)=Id_{M},
\]

\[
z_{0}^{-1}\text{\ensuremath{\delta}}\left(\frac{z_{1}-z_{2}}{z_{0}}\right)Y_{M}\left(u,z_{1}\right)Y_{M}\left(v,z_{2}\right)-z_{0}^{-1}\delta\left(\frac{z_{2}-z_{1}}{-z_{0}}\right)Y_{M}\left(v,z_{2}\right)Y_{M}\left(u,z_{1}\right)
\]

\begin{equation}
=z_{2}^{-1}\left(\frac{z_{1}-z_{0}}{z_{2}}\right)^{-r/T}\delta\left(\frac{z_{1}-z_{0}}{z_{2}}\right)Y_{M}\left(Y\left(u,z_{0}\right)v,z_{2}\right),\label{Jacobi for twisted V-module}
\end{equation}
where $\delta\left(z\right)=\sum_{n\in\mathbb{Z}}z^{n}$ and all binomial expressions (here and below) are to be expanded in nonnegative integral powers of the second variable.
\end{definition}

We use $\Z_{+}$ to denote the set of nonnegative integers.
\begin{definition}
A $g$-\emph{twisted $V$-module} is a weak $g$-twisted $V$-module
$M$ which carries a $\mathbb{C}$-grading induced by the spectrum
of $L(0)$ where $L(0)$ is the component operator of $Y(\omega,z)=\sum_{n\in\mathbb{Z}}L(n)z^{-n-2}.$
That is, we have $M=\bigoplus_{\lambda\in\mathbb{C}}M_{\lambda},$
where $M_{\lambda}=\left\{ w\in M|L(0)w=\lambda w\right\} $. Moreover,
$\dim M_{\lambda}$ is finite and for fixed $\lambda,$ $M_{\frac{n}{T}+\lambda}=0$
for all small enough integers $n.$
\end{definition}

\begin{definition}
An \emph{admissible $g$-twisted $V$-module} $M=\oplus_{n\in\frac{1}{T}\mathbb{Z}_{+}}M\left(n\right)$
is a $\frac{1}{T}\mathbb{Z}_{+}$-graded weak $g$-twisted module
such that $u_{m}M\left(n\right)\subseteq M\left(\mbox{wt}u-m-1+n\right)$
for homogeneous $u\in V$ and $m,n\in\frac{1}{T}\mathbb{Z}.$

\end{definition}
If $g=Id_{V}$ we have the notions of weak, ordinary and admissible
$V$-modules \cite{DLM1}.  From now on, if we say $M$ is a $V$-module, we mean $M$ is an
ordinary $V$-module.
\begin{definition}
A vertex operator algebra $V$ is called \emph{$g$-rational}
if the admissible $g$-twisted module category is semisimple. $V$
is called \emph{rational} if $V$ is $1$-rational.
\end{definition}

We say $V$ is \emph{simple} if as a $V$-module, $V$ is irreducible.
If $V$ is both rational and CFT type, then it is easy to see that $V$ is simple.

\begin{definition}
A vertex operator algebra $V$ is called \emph{$C_{2}$-cofinite} if $\dim V/C_{2}(V)<\infty$, where $C_{2}(V)=\langle u_{-2}v | u, v \in V\rangle$.
\end{definition}
It is proved in \cite{DLM1} that if $V$ is rational, then there are finitely many inequivalent irreducible admissible modules $M^{0}, M^{1},\ldots, M^{d}$
and each irreducible admissible module is an ordinary module. Each $M^{i}$ has weight space decomposition
\begin{align*}
M^{i}=\oplus_{n\in\Z_{+}}M^{i}_{\lambda_{i}+n},
\end{align*}
where $\lambda_{i}\in \C$ is a complex number such that $M^{i}_{\lambda_{i}}\not=0$ and $M^{i}_{\lambda_{i}+n}$ is the eigenspace of $L(0)$ with
eigenvalue $\lambda_{i}+n$. The $\lambda_{i}$ is called the \emph{conformal weight} of $M^{i}$. If $V$ is both rational and $C_{2}$-cofinite, then each $\lambda_{i}$ and the central charge of $V$ are rational numbers \cite{DLM2}.

\begin{definition}
 A vertex operator algebra $V$ is said to be \emph{regular} if any weak $V$-module $M$ is a direct sum of irreducible  $V$-modules.
\end{definition}
 It is shown in \cite{KL} and \cite{ABD} that the regularity is equivalent to rationality and $C_2$-cofiniteness.

\subsection{Invariant bilinear form}
Let $M=\oplus_{\lambda \in\mathbb{C}}M_{\lambda}$ be a $V$-module.
The restricted dual of $M$ is defined by $M'=\oplus_{\lambda \in\mathbb{C}}M_{\lambda}^{\ast}$
where $M_{\lambda}^{\ast}=\text{Hom}_{\mathbb{C}}\left(M_{\lambda},\mathbb{C}\right).$
It is proved in \cite{FHL} that $M'=\left(M',Y_{M'}\right)$ is
naturally a $V$-module such that
\[
\left\langle Y_{M'}\left(v,z\right)f,u\right\rangle =\left\langle f,Y_{M}\left(e^{zL\left(1\right)}\left(-z^{-2}\right)^{L\left(0\right)}v,z^{-1}\right)u\right\rangle ,
\]
for $v\in V,$ $f\in M'$ and $u\in M$, and $\left(M'\right)'\cong M$.
Moreover, if $M$ is irreducible, so is $M'$. A $V$-module $M$
is said to be \emph{self-dual} if $M\cong M'$.
\begin{definition}
A bilinear form $\langle\, ,\, \rangle$ on a $V$-module $M$ is said to be \emph{invariant} if it satisfies the condition
\begin{align*}
\langle Y(a,z)u, v\rangle=\langle u, Y(e^{zL(1)}(-z^{-2})^{L(0)}a, z^{-1})v\rangle
\end{align*}
for $a\in V, u,v \in M$.
\end{definition}
The following result about invariant bilinear forms on $V$ is from \cite{L}:
\begin{theorem}\label{bilinear form}The space of invariant bilinear
forms on $V$ is isomorphic to the space
\[
\left(V_{0}/L\left(1\right)V_{1}\right)^{*}=\mbox{Hom}_{\mathbb{C}}\left(V_{0}/L\left(1\right)V_{1},\mathbb{C}\right).
\]
In particular, if $V$ is a simple vertex operator algebra of CFT type with $V_{1}=0$, then there is a unique nondegenerate invariant bilinear form $\langle\, ,\, \rangle$ on $V$ satisfying $\langle \mathbf{1},\mathbf{1}\rangle=1$.
 \end{theorem}

\subsection{Intertwining operators and fusion rules}

\begin{definition} Let $\left(V,Y\right)$ be a vertex operator algebra
and let $\left(M^{i},Y^{i}\right),\ \left(M^{j},Y^{j}\right)$ and
$\left(M^{k},Y^{k}\right)$ be three $V$-modules. An \emph{intertwining
operator of type $\left(\begin{array}{c}
M^{k}\\
M^{i}\ M^{j}
\end{array}\right)$} is a linear map
\begin{gather*}
\mathcal{Y}\left(\cdot,z\right):\ M^{i}\to\text{\ensuremath{\mbox{Hom}\left(M^{j},\ M^{k}\right)\left\{ z\right\} }}\\
u\mapsto\mathcal{Y}\left(u,z\right)=\sum_{n\in\mathbb{Q}}u_{n}z^{-n-1}
\end{gather*}
satisfying:

(1) For any $u\in M^{i}$ and $v\in M^{j}$, $u_{n}v=0$ for $n$
sufficiently large;

(2) $\mathcal{Y}(L(-1)v,\ z)=\left(\frac{d}{dz}\right)\mathcal{Y}\left(v,z\right)$
for $v\in M^{i}$;

(3) (Jacobi Identity) For any $u\in V,\ v\in M^{i}$,
\begin{alignat*}{1}
 & z_{0}^{-1}\delta\left(\frac{z_{1}-z_{2}}{z_{0}}\right)Y^{k}\left(u,z_{1}\right)\mathcal{Y}\left(v,z_{2}\right)-z_{0}^{-1}\delta\left(\frac{-z_{2}+z_{1}}{z_{0}}\right)\mathcal{Y}\left(v,z_{2}\right)Y^{j}\left(u,z_{1}\right)\\
 & =z_{2}^{-1}\left(\frac{z_{1}-z_{0}}{z_{2}}\right)\mathcal{Y}\left(Y^{i}\left(u,z_{0}\right)v,z_{2}\right).
\end{alignat*}
The space of all intertwining operators of type $\left(\begin{array}{c}
M^{k}\\
M^{i}\ M^{j}
\end{array}\right)$ is denoted by $I_{V}\left(\begin{array}{c}
M^{k}\\
M^{i}\ M^{j}
\end{array}\right)$. Without confusion, we also denote it by $I_{i,j}^{k}.$ Let $N_{i,\ j}^{k}=\dim I_{i,j}^{k}$.
These integers $N_{i,j}^{k}$ are called the \emph{fusion rules}.
\end{definition}

Let $V^{1}$ and $V^{2}$ be vertex operator algebras. Let \{$M^{i}$
, $i=1,2,3$\} be $V^{1}$-modules, and \{$N^{i}$, $i=1,2,3$\} be $V^{2}$-modules.
Then \{$M^{i}\otimes N^{i}$, $i=1,2,3$\} are $V^{1}\otimes V^{2}$-modules
by \cite{FHL}. The following property was given in \cite{ADL}:

\begin{proposition} \label{fusion of tensor product}If $N_{M^{1},M^{2}}^{M^{3}}<\infty$
or $N_{N^{1},N^{2}}^{N^{3}}<\infty,$ then
\[
N_{M^{1}\otimes N^{1},M^{2}\otimes N^{2}}^{M^{3}\otimes N^{3}}=N_{M^{1},M^{2}}^{M^{3}}N_{N^{1},N^{2}}^{N^{3}}.
\]
\end{proposition}

\begin{definition}
Let $M^1$ and $M^2$ be $V$-modules. A \emph{fusion product} for the
ordered pair $\left(M^{1},M^{2}\right)$ is a pair $\left(M,\mathcal{Y}\left(\cdot,z\right)\right)$
which consists of a $V$-module $M$ and an intertwining operator
$\mathcal{Y}\left(\cdot,z\right)$ of type $\left(\begin{array}{c}
W\\
M^{1}\ M^{2}
\end{array}\right)$
satisfies the following universal property: For any $V$-module
$W$ and any intertwining operator $\mathcal{I}\left(\cdot,z\right)$ of type
$\left(\begin{array}{c}
W\\
M^{1}\ M^{2}
\end{array}\right)$, there exists a unique $V$-homomorphism
$\phi$ from $M$ to $W$ such that $\mathcal{I}\left(\cdot,z\right)=\phi\circ \mathcal{Y}\left(\cdot,z\right).$
\end{definition}
From the definition it is easy to see that if a fusion product of $M^{1}$
and $M^{2}$ exists, it is unique up to isomorphism. In this case,
we denote the  fusion product by $M^{1}\boxtimes_{V}M^{2}.$

\begin{definition}
Let $V$ be a simple VOA. A simple $V$-module $M$ is called a \emph{simple current} if for any irreducible $V$-module $W$, $M\boxtimes_V W$ exists and is also a simple $V$-module.
\end{definition}

\subsection{Quantum dimensions and quantum Galois theory}

Now we recall quantum Galois theory and quantum dimensions from \cite{DM2}, \cite{DLM3} and \cite{DJX}.

For a subgroup $G\leq \mathrm{Aut}(V),$ the fixed point set $V^{G}=\{v\in V| g(v)=v, \mathrm{for\, any\,} g \in G\}$ has a vertex operator algebra structure. By \cite{DM2} \cite{DLM3}, we have the following:
\begin{theorem} \label{classical galois theory} Suppose that $V$
is a simple vertex operator algebra and that $G$ is a finite group of automorphisms of $V$. Then the following
hold:

(i) $V=\oplus_{\chi\in\text{Irr}\left(G\right)}V^{\chi},$ where $V^{\chi}$ is the subspace of $V$ on which $G$ acts according to
the character $\chi$. Each $V^{\chi}$ is nonzero;

(ii) For $\chi\in\text{Irr}\left(G\right)$, each $V^{\chi}$ is a
simple module for the $G$-graded vertex operator algebra $\mathbb{C}G\otimes V^{G}$
of the form

\[
V^{\chi}=M_{\chi}\otimes V_{\chi},
\]
where $M_{\chi}$ is the simple $G$-module affording $\chi$
and where $V_{\chi}$ is a simple $V^{G}$-module.

(iii) The map $M_{\chi}\mapsto V_{\chi}$ is a bijection from the
set of inequivalent simple $G$-modules to the set of inequivalent
simple $V^{G}$-modules which are contained in $V$.

\end{theorem}

Now we recall the notion of quantum dimensions from \cite{DJX}. Let $M=\oplus_{n\in \Z_{+}}M_{\lambda+n}$ be a $V$-module. The \emph{formal character} of $M$ is defined to be
\[
\mathrm{ch}_{q}M=\mathrm{tr}_{M} q^{L(0)-c/24}= q^{\lambda-c/24}\sum\limits_{n\in \Z_{+}}(\dim{M_{\lambda+n}})q^{n}.
\]
It is proved in \cite{Z} and \cite{DLM2} that $\mathrm{ch}_{q}M$ converges to a holomorphic function on the domain $|q|<1$ if $V$ is $C_{2}$-cofinite.
We sometimes also use $Z_{M}(\tau)$ to denote the holomorphic function $\mathrm{ch}_{q}M$ with variable $\tau$  in the complex upper half-plane
$\mathbb{H}$ and $q=e^{2\pi i\tau}$. By \cite{DJX}, we have the following:
\begin{definition}
Let $M$ be a $V$-module such that $Z_{V}(\tau)$ and $Z_{M}(\tau)$ exist. The \emph{quantum dimension} of $M$ over $V$ is defined as
\[
q\dim_{V}M=\lim_{y\to 0}\frac{Z_{M}(iy)}{Z_{V}(iy)},
\]
where $y$ is real and positive. Sometimes we use an alternative definition which involves the $q$-characters:
\[
q\dim_{V}M=\lim_{q\to 1^{-}}\frac{\mathrm{ch}_{q}M}{\mathrm{ch}_{q}V}.
\]
\end{definition}

The following results are obtained in \cite{DJX}, \cite{ADJR}:
\begin{proposition} \label{qdim of simple current}
Let $V$ be a
rational and $C_{2}$-cofinite simple vertex operator algebra of CFT
type with $V\cong V'$. Let $M^{0},\ M^{1},\ \cdots,\ M^{d}$ be all
the inequivalent irreducible $V$-modules with $M^{0}\cong V$. Suppose the conformal weights of $M^i$ $(1\leq i\leq d)$ are positive.
Then a $V$-module $M^i$  is a simple current if and only if $q\dim_{V}M^i=1$.
\end{proposition}

\begin{theorem}\label{quantum dimension and orbifold module} Let
$V$ be a rational and $C_{2}$-cofinite simple vertex operator algebra.
Assume $V$ is $g$-rational and the conformal weight of any irreducible
$g$-twisted $V$-module is positive except for $V$ itself for all
$g\in G$. Then
\[
q\dim_{V^{G}}V_{\chi}=\dim W_{\chi}.
\]
\end{theorem}

\begin{remark} \label{product property of qdim}
Let $U$
and $V$ be  vertex operator algebras under the same assumption of Proposition \ref{qdim of simple current}, $M$ be a $U$-module
and $N$ be a $V$-module. Then

\[
q\dim_{U\otimes V}M\otimes N=q\dim_{U}M\cdot q\dim_{V}N.
\]
\end{remark}

\subsection{\label{subsec:The-unitary-series}The unitary series of the Virasoro
VOAs}

 Let
\begin{alignat*}{1}
& c_m:=1-\frac{6}{(m+2)(m+3)}, \ \ \ m=1,2,\dots,\\
& h^{(m)}_{r,s}:=\frac{[r(m+3)-s(m+2)]^2-1}{4(m+2)(m+3)},\ \ \ 1\leq s\leq r\leq m+1.
\end{alignat*}
$L(c_m,h^{(m)}_{r,s})$ is
the irreducible highest weight representation of the Virasoro algebra $L$ with highest weight $\left(c_{m},h^{(m)}_{r,s}\right).$ It is shown in \cite{W} that $L(c_m,0)$ is a rational VOA and $L(c_m,h^{(m)}_{r,s}),$ $ 1\leq s\leq r\leq m+1,$ are all irreducible $L(c_m,0)$-modules. This is the so-called unitary series of the Virasoro VOAs.

The fusion rules among $L(c_m,0)$-modules are computed in \cite{W} and given by
\begin{alignat}{1}\label{fusion rules of virasoro modules}
L(c_m,h^{(m)}_{r_1,s_1})\boxtimes L(c_m,h^{(m)}_{r_2,s_2})=\sum_{i\in I, j\in J} L(c_m,h^{(m)}_{|r_1-r_2|+2i-1,|s_1-s_2|+2j-1}),
\end{alignat}
where
\begin{alignat*}{1}
& I=\{1,2,\dots, \min \{r_1,r_2,m+2-r_1,m+2-r_2\}\}\\
& J=\{1,2,\dots, \min \{s_1,s_2,m+3-s_1,m+3-s_2\}\}.
\end{alignat*}

\begin{definition} Let $V$ be a VOA. A vector $e\in V_{2}$ is called a \emph{conformal
vector with the central charge $c_{e}$ }if it satisfies\emph{ $e_{1}e=2e$
}and $e_{3}e=\frac{c_{e}}{2}\mathbf{1}$. Then the operators $L_{n}^{e}:=e_{n+1},\ n\in\mathbb{Z}$,
satisfy the Virasoro commutation relation
\[
\left[L_{m}^{e},\ L_{n}^{e}\right]=\left(m-n\right)L_{m+n}^{e}+\delta_{m+n,\ 0}\frac{m^{3}-m}{12}c_{e}
\]
for $m,\ n\in\mathbb{Z}.$ A conformal vector $e\in V_{2}$ with the central
charge $\frac{1}{2}$ is called an \emph{Ising vector }if $e$ generates the simple
Virasoro vertex operator algebra $L(\frac{1}{2}, 0)$.
\end{definition}

The fusion rules among $L(c_m,0)$-modules give rise to an involutive automorphism of a VOA.
\begin{theorem}[\cite{M1}]\label{tau involution}
Let $V$ be a VOA and $e\in V$ be a simple Virasoro vector with a central charge $c_m$. Denote by $V_e [h^{(m)}_{r,s}]$ the sum of irreducible $\mathrm{VA}(e) \cong L(c_m,0)$-submodules isomorphic to $L(c_m, h^{(m)}_{r,s}), 1\leq s \leq r \leq m + 1$. Then the linear map
\[
\tau_{e}=
\begin{cases}
(-1)^{r+1} &\hbox{ on } V_e [h^{(m)}_{r,s}]\  \hbox{ if } m \hbox{ is even},  \\
(-1)^{s+1} &\hbox{ on } V_e [h^{(m)}_{r,s}]\  \hbox{ if } m \hbox{ is odd},
\end{cases}
\]
defines an automorphism of $V$ called the $\tau$-involution associated to $e$.
\end{theorem}

We will later consider conformal vectors with central charge $1/2$ and $7/10.$ So we give the $\tau$-involution associated to $c=1/2$ and $c=7/10$ conformal elements precisely. For an Ising vector $e\in V$, one has $$V=V_e[0]\oplus V_e[1/2]\oplus V_e[1/16].$$
By Theorem \ref{tau involution}, the linear map $\tau_{e}$ on $V$ is defined in the following way:
\[
\tau_{e}(v)=
\begin{cases}
1& v\in V_{e}[0]\oplus V_{e}[1/2], \\
-1& v \in V_{e}[1/16].
\end{cases}
\]
For a conformal vector $u\in V$ with a central charge $7/10$, we have
$$V=V_u[0]\oplus V_u[3/2]\oplus V_u[1/10]\oplus V_u[3/5]\oplus V_u[7/16]\oplus V_u[3/80].$$
The $\tau$-involution is defined as
\[
\tau_{u}(v)=
\begin{cases}
1& v\in V_u[0]\oplus V_u[3/2]\oplus V_u[1/10]\oplus V_u[3/5], \\
-1& v \in V_u[7/16]\oplus V_u[3/80].
\end{cases}
\]

 For later use, we recall the definition of $\sigma$-type $c=1/2$ and $c=7/10$ conformal vectors. The corresponding $\sigma$-involutions will also be defined.
\begin{definition}
 An Ising vector $e$ of a VOA $V$ is said to be \emph{of $\sigma$-type} on $V$ if $\tau_e=\hbox{id}$ on $V$.
\end{definition}

In this case, one has $V=V_e[0]\oplus V_e[1/2]$ and the map $\sigma_e$ defined by
\[
\sigma_{e}:=
\begin{cases}
1  &\hbox{ on } V_e [0],  \\
-1  &\hbox{ on } V_e [1/2]
\end{cases}
\]
is an automorphism of $V$ \cite{M1}.

\begin{definition}
A conformal vector $u$ of a VOA $V$ with central charge $7/10$ is said to be \emph{of $\sigma$-type} on $V$ if $V_u[7/16]=V_u[3/80]=0.$
\end{definition}
Let $u\in V$ be a conformal vector of $\sigma$-type with $c=7/10$. Then one has
$$V=V_u[0]\oplus V_u[3/2]\oplus V_u[1/10]\oplus V_u[3/5].$$
Define
\[
\sigma_{u}:=
\begin{cases}
1 &\hbox{ on } V_u [0]\oplus V_u[3/5],  \\
-1 &\hbox{ on } V_u [3/2]\oplus V_u[1/10].
\end{cases}
\]
The fusion rules imply that $\sigma_u$ is an automorphism of $V.$

The quantum dimensions of $L(c_m, h^{(m)}_{r,s})$ are given in \cite{DJX}:
\begin{alignat}{1}\label{qdimVirasoro}
q\dim_{L\left(c_{m},0\right)}{L\left(c_{m},h^{(m)}_{r,s}\right)}=(-1)^{r+s}
\frac{\sin{\frac{\pi r(m+3)}{m+2}}\sin{\frac{\pi s(m+2)}{m+3}}}{\sin{\frac{\pi (m+3)}{m+2}}\sin{\frac{\pi (m+2)}{m+3}}}.
\end{alignat}


\subsection{Griess algebras generated by two Ising vectors}

Next we discuss some results about the weight two space $V_{2}$, where we assume $V$ is a VOA over the real number field $\mathbb{R}$ which is CFT and moonshine type. Then by \cite{L}, there is a unique symmetric invariant bilinear form $\langle \ , \ \rangle$ on V such that $\langle \textbf{1},\textbf{1} \rangle=1$. We also assume that  $\langle \ , \ \rangle$ is positive definite.
 For any two elements $x, y \in V_{2}$, if we define the product $xy:=x_{1}y$, then $V_{2}$ becomes a commutative nonassociative algebra, which is called the \emph{Griess algebra}  \cite{G}. Besides, $V_{2}$ has a bilinear form which is the restriction of $\langle \ , \ \rangle$ on $V_{2}$. Furthermore, for any $x, y, z \in V_{2}$, we have $\langle x, y\rangle\textbf{1}=x_{3}y$, $\langle xy, z \rangle=\langle y, xz\rangle$.
\begin{remark}\label{remarkfixed}
From the relation $\langle xy, z \rangle=\langle y, xz\rangle$, we can easily deduce that for any Ising vector $e\in V$, $\langle \tau_{e}(x), \tau_{e}(y)\rangle=\langle x, y \rangle$ for any $x, y \in V_{2}$.
\end{remark}

The following lemma will be needed later.

\begin{lemma}{\color{blue} }\cite{M2}\label{miyamotodecomposition}
Let $e$ be an Ising vector. Then $V_{2}$ decomposes into
\[
V_{2}=\mathbb{R}e\oplus E^{e}(0)\oplus E^{e}(\frac{1}{2})\oplus E^{e}(\frac{1}{16}),
\]
where $E^{e}(h)$ denotes the eigenspace of $e_{1}$ with eigenvalue $h$.
\end{lemma}

Now let $e, f$ be two Ising vectors in $V_{2}$. We use $\mathrm{VA}(e,f)$ to denote the subVOA generated by the two Ising vectors $e$ and $f$ over $\R$ and let $\C \mathrm{VA}(e,f)$ be its complexification. Let $\G\subseteq \mathrm{VA}(e,f)_{2}$ be the Griess subalgebra generated by $e$ and $f$.  For any automorphsim $\sigma$ of $\mathrm{VA}(e,f)$, we use $e^{\sigma}$ to denote the action of $\sigma$ on $e$. For any two elements $a, y \in \G$, define $\alpha(a, y):=ay-\frac{1}{16}(a+y)$. The following results were given in \cite{S}:
\begin{lemma}\label{sakumaspan}
$\G$ is spanned by
\[
S:=\{e, e^{\tau_{f}}, e^{\tau_{f}\tau_{e}}, f, f^{\tau_{e}}, f^{\tau_{e}\tau_{f}}, \alpha(e, f), \alpha(e, e^{\tau_{f}})\}.
\]
\end{lemma}

\begin{remark}\label{remarkfixed2}
$\alpha(e, f)$ and $\alpha(e, e^{\tau_f})$ are fixed by both $\tau_{e}$ and $\tau_{f}$.
\end{remark}

Let $T$ be the subgroup of $\mathrm{Aut}(V)$ generated by $\tau_{e}$ and $\tau_{f}$, $x^{T}$ denotes the orbit of $x \in V_{2}$ under the action of $T$. Let $\rho=\tau_{e}\tau_{f}$. Then:
\begin{lemma}\cite{S}\label{Sakumalemma}\\
(1) $|e^{T}|=|f^{T}|$. In particular, $e=e^{\rho^n}$ if and only if $f=f^{\rho^n}$.\\
(2) $e^{T}=f^{T}$ if and only if $|e^{T}|$ is odd and $f=e^{\rho^{\frac{n+1}{2}}}$, where $n=|e^{T}|$.\\
(3) $(\tau_{e}\tau_{f})^{|e^{T}\bigcup f^{T}|}=1$ as an automorphism of $V$.
\end{lemma}

\begin{theorem}\cite{S}\label{Sakumathm}
Let $N=|e^{T}\bigcup f^{T}|$.\\
(1) If $N=2$, then $\langle e, f\rangle=0$ or $\frac{1}{2^5}$.\\
(2) If $N=3$, then $\langle e, f\rangle=\frac{13}{2^{10}}$ or $\frac{1}{2^8}$.\\
(3) If $N=4$, then $(\langle e, f\rangle, \langle e, e^{\tau_{f}}\rangle)=(\frac{1}{2^7}, 0)$ or $(\frac{1}{2^8}, \frac{1}{2^5})$.\\
(4) If $N=5$, then $\langle e, f\rangle=\langle e, e^{\tau_{f}}\rangle=\frac{3}{2^9}$.\\
(5) If $N=6$, then $\langle e, f\rangle=\frac{5}{2^{10}}$, $\langle e, e^{\tau_{f}}\rangle=\frac{13}{2^{10}}$ and $\langle e^{\tau_{f}}, f^{\tau_{e}}\rangle=\frac{1}{2^5}$.
\end{theorem}

\begin{remark}
For each case in Theorem \ref{Sakumathm}, the structure of $\G$ can be fully determined by direct calculation following the formulae given in \cite{S}, which is needed in the next section. We recall those formulae from \cite{S}.
\end{remark}


\begin{proposition}\label{S-computation}
Let $e$ and $f$ be two Ising vectors in $V_{2}$, and $\rho=\tau_{e}\tau_{f}.$  Then:\\
(1) For any $n\in \Z$, we have $\langle e, e^{\rho^{n}} \rangle=\langle f, f^{\rho^{n}} \rangle.$\\
(2) For $a\in \{e, f\}, y\in \mathrm{VA}(e, f)_{2}$, we have
\begin{alignat*}{1}
& \langle a, \alpha(a, y) \rangle 
= \frac{31}{16}\langle a, y \rangle-\frac{1}{2^6}, \nonumber\\
& a\cdot\alpha(a, y)=\frac{7}{16}\alpha(a, y)+\left(12\langle a, y\rangle-\frac{25}{2^8}\right)a+\frac{7}{2^9}(y+y^{\tau_{a}}).\nonumber
\end{alignat*}
(3)
We also have the following:
\begin{alignat*}{1}
&\frac{1}{2^{4}}(e-f)- \frac{1}{2^{4}}(f^{\tau_{e}}-e^{\tau_{f}})+\frac{1}{2^4}(e^{\tau_{f}\tau_{e}}-f^{\tau_{e}\tau_{f}})-(\alpha(e, e^{\tau_{f}})-\alpha(f, f^{\tau_{e}}))=0,\label{proequation2},\\
& f\cdot \alpha(e, e^{\tau_{f}})=-\frac{1}{48}(e+e^{\tau_{f}})-\frac{7}{3\cdot 2^{8}}e^{\tau_{f}\tau_{e}}-\frac{13}{2^{8}}f+
\frac{7}{2^{8}}(f^{\tau_{e}}+f^{\tau_{e}\tau_{f}})-\frac{3}{8}\alpha(e, f)+\frac{7}{48}\alpha(e, e^{\tau_{f}}),\\
& e\cdot \alpha(f, f^{\tau_{e}})=-\frac{1}{48}(f+f^{\tau_{e}})-\frac{7}{3\cdot 2^{8}}f^{\tau_{e}\tau_{f}}-\frac{13}{2^{8}}e+
\frac{7}{2^{8}}(e^{\tau_{f}}+e^{\tau_{f}\tau_{e}})-\frac{3}{8}\alpha(e, f)+\frac{7}{48}\alpha(f, f^{\tau_{e}}),\\
& \alpha(e,f)\cdot  \alpha(e,f)  = \frac{7}{3\cdot 2^{11}}(e+e^{\tau_{f}}+e^{\tau_{f}\tau_{e}})+\frac{7}{2^{13}}(f+f^{\tau_{e}}+f^{\tau_{e}\tau_{f}})-\frac{17}{2^{8}} \alpha(e,f)-\frac{7}{3\cdot 2^{9}}\alpha(e, e^{\tau_{f}}), \\
& \alpha(e,f)\cdot \alpha(e, e^{\tau_{f}})  =-\frac{35}{2^{13}}(e+e^{\tau_{f}}+e^{\tau_{f}\tau_{e}})+\frac{7}{2^{12}}(f+f^{\tau_{e}}+f^{\tau_{e}\tau_{f}})-\frac{21}{2^{8}} \alpha(e,f)+\frac{15}{2^{9}}\alpha(e, e^{\tau_{f}}), \\
& \alpha(e, e^{\tau_{f}})\cdot \alpha(e, e^{\tau_{f}}) = \frac{147}{2^{13}}(e+e^{\tau_{f}}+e^{\tau_{f}\tau_{e}})-\frac{63}{2^{9}}\alpha(e, e^{\tau_{f}}).
\end{alignat*}

\end{proposition}

\section{The VOA generated by two Ising vectors $e$ and $f$}

In this section, we will establish the uniqueness of the vertex operator algebras generated by two Ising vectors $e$ and $f$ with $\langle e, f\rangle=\frac{5}{2^{10}}.$ We prove such a vertex operator algebra is isomorphic to the $6A$-algebra $\mathcal{U}_{6A}$ constructed in \cite{LYY}. We also discuss the uniqueness for the cases when $\langle e, f\rangle=\frac{1}{4}, \frac{1}{2^{5}}, 0 $ and the case $\langle e, f\rangle=\frac{1}{2^{8}}$ with  $(\tau_{e}\tau_{f})^4=1$.

Throughout this section, we will treat $V=(V,Y,\textbf{1},\omega)$ be a VOA over the real number field $\mathbb{R}$ which is CFT and moonshine type. Then by \cite{L}, there is a unique symmetric invariant bilinear form $\langle \ , \ \rangle$ on V such that $\langle \textbf{1},\textbf{1} \rangle=1$. We also assume that  $\langle \ , \ \rangle$ is positive definite.
\begin{remark}\label{simpleremark}
Here the reason why we assume $\langle \ , \ \rangle$ is positive definite on $V$ is to guarantee $\langle \ , \ \rangle$ is nondegenerate on any subVOA of $V$.
Then by \cite{L}, any subVOA of $V$ is simple.
\end{remark}
For a VOA $V$ over $\R$, we use $\C V$ to denote its complexification $\C \otimes_{\R} V$. The following results show that there there is no essential difference between the representation theory of
$V$ over $\R$ and that of $\C V$ over $\C$ (see \cite{M}):
\begin{theorem}
Let $V $ be a VOA over $\R$ and let $U$ be an irreducible $\C V$-module with real weights; then $U$ is an irreducible $V$-module or there is a unique
$V$-module $U$ such that $\C U\cong U$ as $\C V$-modules.
\end{theorem}
\begin{corollary}
Let $c \in \R$ and $L(c,0)_{\R}$, $L(c,0)$ be the simple Virasoro vertex operator algebra over $\R$, $\C$ respectively. Then $\C L(c,0)_{\R}$ is isomorphic to $L(c,0)$ as a VOA. Furthermore,  assume that $L(c, h)$ is an irreducible $L(c, 0)$-module with highest weight $h$. If $h\in \R$, then there is a unique irreducible $L(c,0)_{\R}$-module $L(c,h)_{\R}$ such that $L(c, h)\cong \C  L(c,h)_{\R}$. If $W$ is an irreducible $L(c,0)_{\R}$-module, then there is an $h\in \R$ such that $W\cong L(c,h)_{\R}$ and $L(c,h)$ is an irreducible $ L(c,0)$-module with
highest weight $h$.
\end{corollary}
\begin{theorem}
If $\C\otimes_{\R} V$ is rational, then so is $V$.
\end{theorem}

 We now  prove the uniqueness of VOAs generated by two Ising vectors case by case.

\subsection{The case $\langle e, f \rangle=\frac{5}{2^{10}}$}\label{6Acasecase}

Now let $e, f$ be two Ising vectors in a VOA $V$ with $\langle e, f \rangle=\frac{5}{2^{10}}$. We use $\mathrm{VA}(e,f)$ to denote the subVOA generated by $e$ and $f$ over $\R$ and let $\C \mathrm{VA}(e,f)$ be its complexification. Let $\G\subseteq \mathrm{VA}(e,f)_{2}$ be the Griess subalgebra generated by $e$ and $f$.

Set
\begin{align*}
x^{1} &=e,\,\,\, x^{2}=e^{\tau_{f}},\,\,\, x^{3}=e^{\tau_{f}\tau_{e}},\\
x^{4} &=f,\,\,\, x^{5}=f^{\tau_{e}},\,\,\, x^{6}=f^{\tau_{e}\tau_{f}},\\
x^{7} &=\alpha(e, f),\,\,\, x^{8}=\alpha(e, e^{\tau_{f}}).
\end{align*}
The elements $x^i, i=1, \cdots, 8,$
 give a spanning set of $\G$ (see Lemma \ref{sakumaspan}).
First we will show that $\{x^i| i=1, \cdots, 8\}$ are linearly independent.
\begin{lemma}\label{A6dimwt2}
The determinant of the matrix $(\langle x^{i}, x^{j}\rangle)_{1\leqslant i, j \leqslant 8}$ is nonzero. Equivalently, $\dim \G=8$
and $\{x^{i}| i=1,\ldots,8\}$ is a basis of $\G$.
\end{lemma}
\begin{proof} The Gram matrix $(\langle x^{i}, x^{j}\rangle)_{1\leqslant i, j \leqslant 8}$ can be fully determined. We will compute each $\langle x^{i}, x^{j}\rangle$. 
It follows from Theorem \ref{Sakumathm} and Lemma \ref{Sakumalemma} that $\langle e, e^{\tau_{f}} \rangle=\frac{13}{2^{10}}$, $\langle e^{\tau_{f}}, f^{\tau_{e}} \rangle=\frac{1}{2^{5}}$ and $\rho^6=1$.
By appropriate use of Proposition \ref{S-computation}, we get:
\begin{equation}
(\langle x^{i}, x^{j}\rangle)_{1\leqslant i, j \leqslant 8}=\left(
\begin{array}{cccccccc}
\frac{1}{4} & \frac{13}{2^{10}} & \frac{13}{2^{10}} & \frac{5}{2^{10}} & \frac{5}{2^{10}} & \frac{1}{2^{5}} & \frac{-101}{2^{14}} & \frac{147}{2^{14}}\\[0.3cm]
\frac{13}{2^{10}} & \frac{1}{4} & \frac{13}{2^{10}} & \frac{5}{2^{10}} & \frac{1}{2^{5}} & \frac{5}{2^{10}} & \frac{-101}{2^{14}} & \frac{147}{2^{14}}\\[0.3cm]
\frac{13}{2^{10}} & \frac{13}{2^{10}} & \frac{1}{4} & \frac{1}{2^{5}} & \frac{5}{2^{10}} & \frac{5}{2^{10}} & \frac{-101}{2^{14}} & \frac{147}{2^{14}}\\[0.3cm]
\frac{5}{2^{10}} & \frac{5}{2^{10}} & \frac{1}{2^{5}} & \frac{1}{4} & \frac{13}{2^{10}} & \frac{13}{2^{10}} & \frac{-101}{2^{14}} & \frac{-93}{2^{14}}\\[0.3cm]
\frac{5}{2^{10}} & \frac{1}{2^{5}} & \frac{5}{2^{10}} & \frac{13}{2^{10}} & \frac{1}{4} & \frac{13}{2^{10}} & \frac{-101}{2^{14}} & \frac{-93}{2^{14}}\\[0.3cm]
\frac{1}{2^{5}} & \frac{5}{2^{10}} & \frac{5}{2^{10}} & \frac{13}{2^{10}} & \frac{13}{2^{10}} & \frac{1}{4} & \frac{-101}{2^{14}} & \frac{-93}{2^{14}}\\[0.3cm]
\frac{-101}{2^{14}} & \frac{-101}{2^{14}} & \frac{-101}{2^{14}} & \frac{-101}{2^{14}} & \frac{-101}{2^{14}} & \frac{-101}{2^{14}} & \frac{773}{2^{19}} & \frac{-1011}{2^{19}}\\[0.3cm]
\frac{147}{2^{14}} & \frac{147}{2^{14}} & \frac{147}{2^{14}} & \frac{-93}{2^{14}} & \frac{-93}{2^{14}} & \frac{-93}{2^{14}} & \frac{-1011}{2^{19}}
& \frac{3717}{2^{19}}\\[0.3cm]
\end{array}
\right).
\end{equation}

The determinant of the matrix is nonzero, thus $\{x^i|i=1,\cdots, 8\}$ forms a basis of $\G .$

We give the calculations for some entries of the above matrix. Similar methods can be applied to other entries. For example,
\begin{alignat}{1}
\langle x^{1}, x^{3} \rangle&=\langle e, e^{\tau_{f}\tau_{e}} \rangle=\langle \tau_{e}(e), \tau_{e}(e^{\tau_{f}\tau_{e}}) \rangle
=\langle e,e^{\tau_{f}} \rangle=\frac{13}{2^{10}},\nonumber\\
\langle x^{4}, x^{5} \rangle&=\langle f, f^{\tau_{e}} \rangle=\langle \tau_{f}(f), \tau_{f}(f^{\tau_{e}}) \rangle
=\langle f,f^{\tau_{e}\tau_{f}} \rangle=\langle x^{1}, x^{3} \rangle=\frac{13}{2^{10}},\nonumber
\end{alignat}

\begin{alignat}{1}
\langle x^4, x^8 \rangle & = \langle f, ee^{\tau_{f}}\rangle-\frac{1}{16}(\langle f, e\rangle+ \langle f, e^{\tau_{f}} \rangle) \nonumber\\
& = \langle ef, e^{\tau_{f}}\rangle-\frac{1}{16}(\langle f, e\rangle+ \langle f, e^{\tau_{f}} \rangle) \nonumber\\
& = \langle \alpha(e, f)+\frac{1}{16}(e+f), e^{\tau_{f}}\rangle-\frac{1}{16}(\langle f, e\rangle+ \langle f, e^{\tau_{f}} \rangle) \nonumber\\
& =  \langle e, \alpha(e, f)\rangle+\frac{1}{16}\langle e, e^{\tau_{f}}\rangle-\frac{1}{16} \langle e, f \rangle, \nonumber\\
\langle x^{7}, x^{8} \rangle& = \langle ef, \alpha(e, e^{\tau_{f}}) \rangle-\frac{1}{16}\langle e, \alpha(e, e^{\tau_{f}})\rangle-\frac{1}{16}\langle f, \alpha(e, e^{\tau_{f}})\rangle \nonumber\\
& = \langle f, e\cdot\alpha(e, e^{\tau_{f}}) \rangle-\frac{1}{16}\langle e, \alpha(e, e^{\tau_{f}})\rangle-\frac{1}{16}\langle f, \alpha(e, e^{\tau_{f}})\rangle \nonumber.
\end{alignat}

\end{proof}

Next we work on the VOA $\C \mathrm{VA}(e, f)$. Since $e$ is an Ising vector, $e^{\tau_{f}}$ is also an Ising vector. Let $\mathcal{V}$ be the subalgebra of $\C\mathrm{VA}(e, f)$
generated by $e$ and $e^{\tau_{f}}$. Since $\langle e, e^{\tau_{f}}\rangle=\frac{13}{2^{10}}$, it follows from \cite{M2}, \cite{SY} that
\begin{alignat}{1}
\mathcal{V}\cong & \left(L\left(\frac{4}{5}, 0\right)\oplus L\left(\frac{4}{5}, 3\right)\right)\otimes  \left(L\left(\frac{6}{7}, 0\right)\oplus L\left(\frac{6}{7}, 5\right)\right) \nonumber\\
& \oplus L\left(\frac{4}{5}, \frac{2}{3}\right)^{+}\otimes L\left(\frac{6}{7}, \frac{4}{3}\right)^{+}\oplus L\left(\frac{4}{5}, \frac{2}{3}\right)^{-}\otimes L\left(\frac{6}{7}, \frac{4}{3}\right)^{-}.
\end{alignat}
\begin{remark}
On $M=L\left(\frac{6}{7}, \frac{4}{3}\right)$, there are two $L\left(\frac{6}{7}, 0\right)\oplus L\left(\frac{6}{7}, 5\right)$-module structures.
Namely, if $Y_{M}(v, z)$ is a vertex operator of $v \in L\left(\frac{6}{7}, 0\right)\oplus L\left(\frac{6}{7}, 5\right)$, then the other is defined by
$Y_{M}(v, z)$ for $v \in L\left(\frac{6}{7}, 0\right)$ and $-Y_{M}(v, z)$ for $v \in L\left(\frac{6}{7}, 5\right)$. We denote them by
$L\left(\frac{6}{7}, \frac{4}{3}\right)^{\pm}$. The definition of $L\left(\frac{4}{5}, \frac{2}{3}\right)^{\pm}$ are similar.
\end{remark}
The following two lemmas are from \cite{SY}:
\begin{lemma}
$\mathcal{V}$ is rational.
\end{lemma}
\begin{lemma}\label{v-module}
All the irreducible $\mathcal{V}$-modules are as follows:
\begin{alignat}{1}
 \mathcal{V}=\mathcal{V}(0)\cong & \left(L\left(\frac{4}{5}, 0\right)\oplus L\left(\frac{4}{5}, 3\right)\right)\otimes  \left(L\left(\frac{6}{7}, 0\right)\oplus L\left(\frac{6}{7}, 5\right)\right) \nonumber\\
& \oplus L\left(\frac{4}{5}, \frac{2}{3}\right)^{+}\otimes L\left(\frac{6}{7}, \frac{4}{3}\right)^{+}\oplus L\left(\frac{4}{5}, \frac{2}{3}\right)^{-}\otimes L\left(\frac{6}{7}, \frac{4}{3}\right)^{-}.\nonumber\\
 \mathcal{V}\left(\frac{1}{7}\right)\cong & \left(L\left(\frac{4}{5}, 0\right)\oplus L\left(\frac{4}{5}, 3\right)\right)\otimes  \left(L\left(\frac{6}{7}, \frac{1}{7}\right)\oplus L\left(\frac{6}{7}, \frac{22}{7}\right)\right) \nonumber\\
& \oplus L\left(\frac{4}{5}, \frac{2}{3}\right)^{+}\otimes L\left(\frac{6}{7}, \frac{10}{21}\right)^{+}\oplus L\left(\frac{4}{5}, \frac{2}{3}\right)^{-}\otimes L\left(\frac{6}{7}, \frac{10}{21}\right)^{-}.\nonumber
\end{alignat}
\begin{alignat}{1}
 \mathcal{V}\left(\frac{5}{7}\right)\cong & \left(L\left(\frac{4}{5}, 0\right)\oplus L\left(\frac{4}{5}, 3\right)\right)\otimes  \left(L\left(\frac{6}{7}, \frac{5}{7}\right)\oplus L\left(\frac{6}{7}, \frac{12}{7}\right)\right) \nonumber\\
& \oplus L\left(\frac{4}{5}, \frac{2}{3}\right)^{+}\otimes L\left(\frac{6}{7}, \frac{1}{21}\right)^{+}\oplus L\left(\frac{4}{5}, \frac{2}{3}\right)^{-}\otimes L\left(\frac{6}{7}, \frac{1}{21}\right)^{-}.\nonumber\\
 \mathcal{V}\left(\frac{2}{5}\right)\cong & \left(L\left(\frac{4}{5}, \frac{2}{5}\right)\oplus L\left(\frac{4}{5}, \frac{7}{5}\right)\right)\otimes  \left(L\left(\frac{6}{7}, 0\right)\oplus L\left(\frac{6}{7}, 5\right)\right) \nonumber\\
& \oplus L\left(\frac{4}{5}, \frac{1}{15}\right)^{+}\otimes L\left(\frac{6}{7}, \frac{4}{3}\right)^{+}\oplus L\left(\frac{4}{5}, \frac{1}{15}\right)^{-}\otimes L\left(\frac{6}{7}, \frac{4}{3}\right)^{-}.\nonumber
\end{alignat}
\begin{alignat}{1}
 \mathcal{V}\left(\frac{19}{35}\right)\cong & \left(L\left(\frac{4}{5}, \frac{2}{5}\right)\oplus L\left(\frac{4}{5}, \frac{7}{5}\right)\right)\otimes  \left(L\left(\frac{6}{7}, \frac{1}{7}\right)\oplus L\left(\frac{6}{7}, \frac{22}{7}\right)\right) \nonumber\\
& \oplus L\left(\frac{4}{5}, \frac{1}{15}\right)^{+}\otimes L\left(\frac{6}{7}, \frac{10}{21}\right)^{+}\oplus L\left(\frac{4}{5}, \frac{1}{15}\right)^{-}\otimes L\left(\frac{6}{7}, \frac{10}{21}\right)^{-}.\nonumber\\
 \mathcal{V}\left(\frac{39}{35}\right)\cong & \left(L\left(\frac{4}{5}, \frac{2}{5}\right)\oplus L\left(\frac{4}{5}, \frac{7}{5}\right)\right)\otimes  \left(L\left(\frac{6}{7}, \frac{5}{7}\right)\oplus L\left(\frac{6}{7}, \frac{12}{7}\right)\right) \nonumber\\
& \oplus L\left(\frac{4}{5}, \frac{1}{15}\right)^{+}\otimes L\left(\frac{6}{7}, \frac{1}{21}\right)^{+}\oplus L\left(\frac{4}{5}, \frac{1}{15}\right)^{-}\otimes L\left(\frac{6}{7}, \frac{1}{21}\right)^{-}.\nonumber
\end{alignat}
\end{lemma}
By \cite{M2}, we have the following proposition:
\begin{proposition} \label{A6prop1} $x^{1}, x^{2}, x^{3}, x^{8}$ are defined as before. Then:\\
(1) $\mathcal{V}_{2}=\mathrm{span}\{x^{1}, x^{2}, x^{3}, x^{8}\}$.\\
(2) Let
\begin{alignat}{1}
\omega^{1} & =\frac{2^6}{135}\left(x^{1}+x^{2}+x^{3}-2^{4}x^{8}\right), \nonumber\\
\omega^{2} & =\frac{2^4}{27}\left(x^{1}+x^{2}+x^{3}+\frac{2^5}{7}x^{8}\right) \nonumber.
\end{alignat}
Then $\omega^{1}$ and $\omega^{2}$ are orthogonal conformal vectors with central charges $\frac{4}{5}, \frac{6}{7}$ respectively.
\end{proposition}
\begin{proposition}\label{A6prop2}
$x^{i}, i=1,\ldots, 8$ are defined as before. Let
\begin{alignat}{1}
\omega^{3}=-\frac{1}{18}\left(x^{1}+x^{2}+x^{3}-21(x^{4}+x^{5}+x^{6})-9\cdot 2^{5}x^{7}-\frac{29\cdot 2^{5}}{7}x^{8}\right).\nonumber
\end{alignat}
Then $\omega^{3}$ is a conformal vector with central charge $\frac{25}{28}$. Furthermore,
$\omega^{1}+\omega^{2}+\omega^{3}$ is the Virasoro vector of $\C \mathrm{VA}(e, f)$ and $x^{i}\omega^{3}=0$ for $i=1, 2, 3, 8$.
\end{proposition}
\begin{proof}
It can be verified directly that $\omega^{3}$ is a conformal vector with central charge $\frac{25}{28},$
$x^{i}\omega^{3}=0$ and $\langle x^i, \omega^3\rangle =0$ for $i=1, 2, 3, 8$ by computing the product $x^{i}x^{j}$ and using the inner product $\langle x^{i}, x^{j}\rangle$ for $1\leq i, j \leq 8$ calculated before.

For the proof of $\omega^{1}+\omega^{2}+\omega^{3}$ is the Virasoro vector of $\C \mathrm{VA}(e, f)$, it suffices to show that $Y((\omega^{1}+\omega^{2}+\omega^{3})(0)v,z)=
\frac{d}{dz}Y(v,z)$ for any $v \in \C \mathrm{VA}(e, f)$. By \cite{M2}, $\omega^{1}+\omega^{2}=e+a$ where $a$ is a highest weight vector of $L(\frac{1}{2}, 0)$. Since $e\omega^{3}=0$, we can write $\omega^{1}+\omega^{2}+\omega^{3}=e+a'$ where $a'=a+\omega^{3}$ is a highest weight vector of $\mathrm{VA}(e)\cong L(\frac{1}{2}, 0)$. So $e(0)a'=0$. By skew-symmetry property, $a'(0)e=0$. So we have
\begin{alignat}{1}
 Y((\omega^{1}+\omega^{2}+\omega^{3})(0)e, z)=Y((e+a')(0)e,z)=Y(e(0)e,z)=
\frac{d}{dz}Y(e,z).\nonumber
\end{alignat}
Similarly, let
\begin{alignat}{1}
\widetilde{\omega^{1}} & =\frac{2^6}{135}\left(x^{4}+x^{5}+x^{6}-2^{4}\alpha(f, f^{\tau_e})\right), \nonumber\\
\widetilde{\omega^{2}} & =\frac{2^4}{27}\left(x^{4}+x^{5}+x^{6}+\frac{2^5}{7}\alpha(f, f^{\tau_e})\right), \nonumber\\
\widetilde{\omega^{3}} &
-\frac{1}{18}\left(x^{4}+x^{5}+x^{6}-21(x^{1}+x^{2}+x^{3})-9\cdot 2^{5}x^{7}-\frac{29\cdot 2^{5}}{7}\alpha(f, f^{\tau_e})\right).\nonumber
\end{alignat}
Then we have
\begin{alignat}{1}
 Y((\widetilde{\omega^{1}}+\widetilde{\omega^{2}}+\widetilde{\omega^{3}})(0)f, z)=\frac{d}{dz}Y(f,z).\nonumber
\end{alignat}
By (4) of Proposition \ref{S-computation},
it is easy to see that $\omega^{1}+\omega^{2}+\omega^{3}=\widetilde{\omega^{1}}+\widetilde{\omega^{2}}+\widetilde{\omega^{3}}$, so we have
\begin{alignat}{1}
 Y((\omega^{1}+\omega^{2}+\omega^{3})(0)f, z)=
\frac{d}{dz}Y(f,z).\nonumber
\end{alignat}
Since $\C\mathrm{VA}(e, f)$ is generated by $e$ and $f$, by \cite{L1}, we have $Y((\omega^{1}+\omega^{2}+\omega^{3})(0)v, z)=
\frac{d}{dz}Y(v,z)$ for any $v \in \C \mathrm{VA}(e, f)$.
\end{proof}

We can find common eigenvectors of $\omega^{1}, \omega^{2},\omega^{3}$ in $\G$. We also compute some products that will be used later.

\begin{lemma}\label{v01234}
Set
\begin{alignat}{1}
 v^{0} & = \frac{16}{27}(-2x^{1}+x^{2}+x^{3}),\nonumber\\
v^{1} & = \frac{16}{9}(x^{2}-x^{3}),\nonumber\\
v^{2} & = \frac{4}{3}\left(-\frac{1}{9}(-2x^{1}+x^{2}+x^{3})+(x^{4}+x^{5}-2x^{6})\right),\nonumber\\
v^{3} & =\frac{4}{9}(x^{2}-x^{3})+4(x^{4}-x^{5}),\nonumber\\
v^{4} & = 13(x^{1}+x^{2}+x^{3})+27(x^{4}+x^{5}+x^{6})+81\cdot 16 \cdot x^{7} +23\cdot 16\cdot x^{8}. \nonumber
\end{alignat}
Then we have the following:
\begin{alignat}{1}
& \omega^{1}v^{0}=\frac{2}{3}v^{0},\,\,\,\,\,\, \omega^{1}v^{1}=\frac{2}{3}v^{1},\,\,\,\,\,\, \omega^{1}v^{2}=\frac{2}{3}v^{2},
\,\,\,\,\,\, \omega^{1}v^{3}=\frac{2}{3}v^{3},\,\,\,\,\,\, \omega^{1}v^{4}=0, \nonumber\\
& \omega^{2}v^{0}=\frac{4}{3}v^{0},\,\,\,\,\,\, \omega^{2}v^{1}=\frac{4}{3}v^{1},\,\,\,\,\,\, \omega^{2}v^{2}=\frac{1}{21}v^{2},
\,\,\,\,\,\, \omega^{2}v^{3}=\frac{1}{21}v^{3},\,\,\,\,\,\, \omega^{2}v^{4}=\frac{5}{7}v^{4}, \nonumber\\
& \omega^{3}v^{0}=0,\,\,\,\,\,\, \omega^{3}v^{1}=0,\,\,\,\,\,\, \omega^{3}v^{2}=\frac{9}{7}v^{2},
\,\,\,\,\,\, \omega^{3}v^{3}=\frac{9}{7}v^{3},\,\,\,\,\,\, \omega^{3}v^{4}=\frac{9}{7}v^{4}, \nonumber\\
&  v^{0}v^{0}=\frac{5}{6}\omega^{1}+\frac{14}{9}\omega^{2}-\frac{10}{9} v^{0},\,\,\,\,\,\,v^{0}v^{1}=\frac{10}{9}v^{1},\,\,\,\,\,\,
v^{0}v^{3}=\frac{1}{9}v^{3}, \,\,\,\,\,\,v^{2}v^{3}=-\frac{5}{9}v^{1}+\frac{22}{9}v^{3},\nonumber\\
& \langle v^{0}, v^{0}\rangle=\frac{1}{2},\,\,\,\,\,\,\langle v^{1}, v^{1}\rangle=\frac{3}{2}
,\,\,\,\,\,\,\langle v^{2}, v^{2}\rangle=\frac{5}{2},\,\,\,\,\,\,\langle v^{3}, v^{3}\rangle=\frac{15}{2}.\nonumber
\end{alignat}
\end{lemma}

For any $L\left(\frac{4}{5}, 0\right)\otimes L\left(\frac{6}{7}, 0\right) \otimes L\left(\frac{25}{28}, 0\right)$-module  $L\left(\frac{4}{5}, h_{1}\right)\otimes L\left(\frac{6}{7}, h_{2}\right) \otimes L\left(\frac{25}{28}, h_{3}\right)$, we simply denote it by $[h_{1}, h_{2}, h_{3}]$. Now we can state our main theorem:
\begin{theorem}\label{A6mainthm}
As a module of $\mathcal{V}\otimes L\left(\frac{25}{28}, 0\right)$,
\begin{alignat}{1}
\C\mathrm{VA}(e, f)\cong \mathcal{V}\otimes L\left(\frac{25}{28}, 0\right)\oplus \mathcal{V}\left(\frac{1}{7}\right)\otimes L\left(\frac{25}{28}, \frac{34}{7}\right)\oplus \mathcal{V}\left(\frac{5}{7}\right)\otimes L\left(\frac{25}{28}, \frac{9}{7}\right).\nonumber
\end{alignat}
\end{theorem}

In order to prove the theorem, one needs to show that as an extension of the VOA
$\mathcal{V}\otimes L\left(\frac{25}{28}, 0\right),$ each irreducible $\mathcal{V}\otimes L\left(\frac{25}{28}, 0\right)$-module in $\C\mathrm{VA}(e, f)$  has multiplicity 1. So we recall some results on extensions of VOAs.

 Let $V^1, V^2$ be  two vertex operator algebras which are regular, self-dual, simple and of CFT-type. Assume $V$ is an extension of $V^1\otimes V^2,$ i.e.
\begin{alignat*}{1}
V=\bigoplus_{i\in I, j \in J}{Z_{ij}M^i\otimes N^j}
\end{alignat*}
where $Z_{ij}(i\in I, j\in J)$ are nonnegative integers which are not all zeros and $\{M^i|i\in I\}$ (respectively $\{N^j|j\in J\}$) are  inequivalent irreducible $V^1$-modules (respectively $V^2$-modules). It is proved \cite{HKL} that $V$ is regular and simple CFT-type.
The following theorem is given in \cite{Lin}:
\begin{theorem}\label{mirror}
Let $V,$ $V^1,$ $V^2$ and $Z_{i,j}$ be as above. Assume that
$\Hom_{V^1\otimes V^2}{(V^1\otimes N^j, V)}= \C$ (respectively
$\Hom_{V^1\otimes V^2}{(M^i\otimes V^2, V)}= \C$) if and only if $N^j=V^2$ (respectively $M^i=V^1$). Then we have
$$Z_{ij}=1 \,\,\, \mathrm{if} \,\,\, Z_{ij}\neq 0.$$
\end{theorem}

Let $\Z_{+}$ be the set of nonnegative integers. Next we will use the following two propositions to establish Theorem \ref{A6mainthm}.
\begin{proposition}\label{A6propm1}
As a $\mathcal{V}\otimes L(\frac{25}{28}, 0)$-module,
\begin{align*}
\C\mathrm{VA}(e,f) \cong \mathcal{V}\otimes L(\frac{25}{28}, 0) \oplus m\cdot\mathcal{V}(\frac{1}{7})\otimes L(\frac{25}{28}, \frac{34}{7}) \oplus \mathcal{V}(\frac{5}{7})\otimes L(\frac{25}{28}, \frac{9}{7}),
\end{align*}
where $m \in \Z_{+}$ is the multiplicity of $\mathcal{V}(\frac{1}{7})\otimes L(\frac{25}{28}, \frac{34}{7})$. In particular, $\C\G=\C\mathrm{VA}(e, f)_{2}$.
\end{proposition}
\begin{proof}
By Proposition \ref{A6prop1} and Proposition \ref{A6prop2}, we have $\mathcal{V}\otimes L(\frac{25}{28}, 0)\subseteq \C\mathrm{VA}(e,f)$.
Note that $\mathcal{V}\otimes L(\frac{25}{28}, 0)$, $\mathcal{V}(\frac{1}{7})\otimes L(\frac{25}{28}, \frac{34}{7})$ and $\mathcal{V}(\frac{5}{7})\otimes L(\frac{25}{28}, \frac{9}{7})$ are the only irreducible $\mathcal{V}\otimes L(\frac{25}{28}, 0)$-modules which have integral weights. So we have
\begin{alignat}{1}\label{weakerv2}
\C\mathrm{VA}(e, f)\cong \mathcal{V}\otimes L\left(\frac{25}{28}, 0\right)\oplus m\cdot\mathcal{V}\left(\frac{1}{7}\right)\otimes L\left(\frac{25}{28}, \frac{34}{7}\right)\oplus n\cdot\mathcal{V}\left(\frac{5}{7}\right)\otimes L\left(\frac{25}{28}, \frac{9}{7}\right)
\end{alignat}
where $m, n \in \Z_+$. By Lemma \ref{v01234}, we have $n\neq 0$. Since both $\mathcal{V}$ and $L\left(\frac{25}{28}, 0\right)$ are regular and simple CFT-type vertex operator algebras, we have $\mathcal{V}\otimes L\left(\frac{25}{28}, 0\right)$ is regular and simple CFT-type vertex operator algebra. By \cite{HKL}, $\C\mathrm{VA}(e, f)$ is also regular and simple CFT-type vertex operator algebra. Since $\C\mathrm{VA}(e, f)$ is a moonshine type vertex operator algebra, by \cite{L}, $\C\mathrm{VA}(e, f)$ is self-dual.
Similarly, $\mathcal{V}$ and $L\left(\frac{25}{28}, 0\right)$ are self-dual.
 By the equation \ref{weakerv2},  $$\Hom_{\mathcal{V}\otimes L\left(\frac{25}{28}, 0\right)}{\left(\mathcal{V}\otimes L\left(\frac{25}{28}, \frac{9}{7}\right), \C\mathrm{VA}(e, f)\right)}=0.$$ Now $n=1$ follows directly from Theorem \ref{mirror}.

Next, we show that $\dim \C\mathrm{VA}(e, f)_{2}=8$, i.e. $\C\G=\C\mathrm{VA}(e, f)_{2}$.
By Lemma \ref{A6dimwt2}, we have $\dim \C\G=8$. It is easy to see that
\begin{alignat}{1}
& \dim \left( \mathcal{V}\otimes L\left(\frac{25}{28}, 0\right)\right)_{2}=5,\nonumber\\
& \dim \left( \mathcal{V}\left(\frac{1}{7}\right)\otimes L\left(\frac{25}{28}, \frac{34}{7}\right)\right)_{2}=0,\nonumber\\
& \dim \left( \mathcal{V}\left(\frac{5}{7}\right)\otimes L\left(\frac{25}{28}, \frac{9}{7}\right)\right)_{2}=3.\nonumber
\end{alignat}
So $\dim \C\mathrm{VA}(e, f)_{2}=8=\dim \C\G$, hence $\C\G=\C\mathrm{VA}(e, f)_{2}$.
\end{proof}

\begin{proposition}\label{A6propm2}
The number $m$ in Proposition \ref{A6propm1} is nonzero. Thus $m=1$.
\end{proposition}
\begin{proof}
The idea of the proof is similar to Theorem 5.3 in \cite{SY}. First we compute the highest weight vector for $[0, 0, 0]$ with highest weight $[3,0,0]$. For any $[0, 0, 0]$-module $U$,
let $U^{(\tau_{e}, \pm)}=\{u\in U| \tau_{e}(u)=\pm u\}$. We claim that $\{\omega^{1}_{0}v^{1},\omega^{2}_{0}v^{1}, v^{0}_{0}v^{1}\}$ is a
set of linearly independent vectors in $\left(\mathcal{V}\otimes L\left(\frac{25}{28}, 0\right)\right)^{(\tau_{e}, -)}_{3}$. Set $t^{1}=\omega^{1}_{0}v^{1}$, $t^{2}=\omega^{2}_{0}v^{1}$ and $t^{3}=v^{0}_{0}v^{1}$. Using the commutator formula $[a_{m}, b_{n}]=\sum_{i \in \Z_{+}}\binom{m}{i}(a_{i}b)_{(m+n-i)}$, an invariant property $\langle a_{m}b^{1},\,  b^{2}\rangle=\langle b^{1}, \, a_{(-m+2)}b^{2}\rangle$ for $a\in \C\mathrm{VA}(e, f)_{2}$,
and an identity $(a_{0}b)_{m}=[a_{1}, b_{(m-1)}]-(a_{1}b)_{(m-1)}$, we can calculate all $\langle t^{i}, \, t^{j}\rangle$, $1\leq i, j\leq 3$. For example,
we compute $\langle t^{1}, \, t^{3}\rangle=\langle \omega^{1}_{0}v^{1}, \, v^{0}_{0}v^{1} \rangle$:
\begin{alignat}{1}
\langle \omega^{1}_{0}v^{1}, \, v^{0}_{0}v^{1} \rangle &= \langle v^{1}, \, \omega^{1}_{2}v^{0}_{0}v^{1} \rangle=
 \langle v^{1}, \, [\omega^{1}_{2}, v^{0}_{0}]v^{1} \rangle \nonumber\\
& = \langle v^{1}, \, \left((\omega^{1}_{0}v^{0})_{2}+2(\omega^{1}_{1}v^{0})_{1}+(\omega^{1}_{2}v^{0})_{0}\right)v^{1} \rangle\nonumber\\
& =\langle v^{1}, \, \left([\omega^{1}_{1}, v^{0}_{1}]+(\omega^{1}_{1}v^{0})_{1}\right)v^{1} \rangle \nonumber\\
& =\langle \omega^{1}_{1}v^{1},\,  v^{0}_{1}v^{1} \rangle-\langle v^{1},\,  v^{0}_{1}\omega^{1}_{1}v^{1} \rangle+ \langle v^{1}, \, (\omega^{1}_{1}v^{0})_{1}v^{1} \rangle \nonumber\\
& =\frac{2}{3}\langle v^{1}, \, v^{0}_{1}v^{1} \rangle =\frac{20}{27}\langle v^{1}, \, v^{1} \rangle=\frac{10}{9}. \nonumber
\end{alignat}
By a similar way, we can compute all $\langle t^{i}, \, t^{j}\rangle$, $1\leq i, j\leq 3$, and it is a routine work to check that $\mathrm{det}\left(\langle t^{i},\, t^{j}\rangle\right)_{1\leq i, j\leq 3}\not=0$. Therefore $t^{1}, t^{2}$ and $t^{3}$ are linearly independent. One can also see that
\begin{alignat}{1}
u^{1}=v^{0}_{0}v^{1}-\frac{5}{9}(\omega^{1}_{0}+\omega^{2}_{0})v^{1}\nonumber
\end{alignat}
is a non-zero highest weight vector for $[0, 0, 0]$ with highest weight $[3, 0, 0]$. Next we compute the highest weight vector for $[0,0,0]$ with highest weight $[0, \frac{12}{7}, \frac{9}{7}]$. We claim $\{\omega^{1}_{0}v^{3},\omega^{2}_{0}v^{3}, \omega^{3}_{0}v^{3}, v^{0}_{0}v^{3}\}$ is a
set of linearly independent vectors in $\left(\mathcal{V}\left(\frac{5}{7}\right)\otimes L\left(\frac{25}{28}, \frac{9}{7}\right)\right)^{(\tau_{e}, -)}_{3}$. Set $r^{1}=\omega^{1}_{0}v^{3}$, $r^{2}=\omega^{2}_{0}v^{3}$, $r^{3}=\omega^{3}_{0}v^{3}$ and $r^{4}=v^{0}_{0}v^{3}$. By a similar method used in computations of $\langle t^{i}, \, t^{j}\rangle$, $1\leq i, j\leq 3$, we can calculate all  $\langle r^{i}, \, r^{j}\rangle$, $1\leq i, j\leq 4$. It is also
 a routine work to check that $\mathrm{det}\left(\langle r^{i},\, r^{j}\rangle\right)_{1\leq i, j\leq 4}\not=0$. Therefore $r^{1}, r^{2}$, $r^{3}$ and $r^{4}$ are linearly independent. One can also see that
\begin{alignat}{1}
u^{2}=v^{0}_{0}v^{3}-\frac{1}{18}\omega^{1}_{0}v^{3}-\frac{14}{9}\omega^{2}_{0}v^{3} \nonumber
\end{alignat}
is a non-zero highest weight vector for $[0, 0, 0]$ with highest weight $[0, \frac{12}{7}, \frac{9}{7}]$.
Next we compute the highest weight vector for $[0,0,0]$ with highest weight $[0, 5, 0]$. Set
\begin{alignat}{1}
& a^{1}=\omega^{1}_{-2}v^{1},\,\,\,\,\, a^{2}=\omega^{1}_{-1}\omega^{1}_{0}v^{1},\,\,\,\,\,
 a^{3}=\omega^{1}_{-1}\omega^{2}_{0}v^{1},\,\,\,\,\, a^{4}=\omega^{1}_{0}\omega^{1}_{0}\omega^{2}_{0}v^{1}, \nonumber \\
& a^{5}=\omega^{1}_{0}\omega^{2}_{-1}v^{1},\,\,\,\,\,a^{6}=\omega^{1}_{0}\omega^{2}_{0}\omega^{2}_{0}v^{1},\,\,\,\,\,
 a^{7}=\omega^{1}_{0}\omega^{3}_{-1}v^{1},\,\,\,\,\, a^{8}=\omega^{2}_{-2}v^{1}, \nonumber \\
& a^{9}=\omega^{2}_{-1}\omega^{2}_{0}v^{1},\,\,\,\,\,a^{10}=\omega^{2}_{0}\omega^{3}_{-1}v^{1},\,\,\,\,\,
 a^{11}=\omega^{3}_{-2}v^{1},\,\,\,\,\, a^{12}=\omega^{1}_{-1}u^{1}, \nonumber\\
& a^{13}=\omega^{1}_{0}\omega^{1}_{0}u^{1},\,\,\,\,\, a^{14}=\omega^{2}_{-1}u^{1},\,\,\,\,\,
a^{15}=\omega^{3}_{-1}u^{1},\,\,\,\,\, a^{16}=v^{0}_{-2}v^{1}.\nonumber
\end{alignat}
By a similar method used in computations of $\langle t^{i}, \, t^{j}\rangle$, $1\leq i, j\leq 3$, we can calculate all  $\langle a^{i}, \, a^{j}\rangle$, $1\leq i, j\leq 16$. It is also
 a routine work to check that $\mathrm{det}\left(\langle a^{i},\, a^{j}\rangle\right)_{1\leq i, j\leq 16}\not=0$. Therefore $\{a^{i}, i=1,\ldots, 16\} $ is a linearly independent set. One can also see that
\begin{alignat}{1}
u^{3}=& \frac{275}{243}\omega^{1}_{-2}v^{1}-\frac{50}{81}\omega^{1}_{-1}\omega^{1}_{0}v^{1}
+\frac{25}{36}\omega^{1}_{-1}\omega^{2}_{0}v^{1}-\frac{25}{72}\omega^{1}_{0}\omega^{1}_{0}\omega^{2}_{0}v^{1}\nonumber\\
& +\frac{14}{45}\omega^{1}_{0}\omega^{2}_{-1}v^{1}-\frac{7}{180}\omega^{1}_{0}\omega^{2}_{0}\omega^{2}_{0}v^{1}
+\frac{140}{243}\omega^{2}_{-2}v^{1}-\frac{7}{81}\omega^{2}_{-1}\omega^{2}_{0}v^{1}\nonumber\\
& +\frac{5}{39}\omega^{1}_{-1}u^{1}+\frac{3}{26}\omega^{1}_{0}\omega^{1}_{0}u^{1}
+\frac{28}{9}\omega^{2}_{-1}u^{1}-v^{0}_{-2}v^{1}\nonumber
\end{alignat}
is a non-zero highest weight vector for $[0, 0, 0]$ with highest weight $[0, 5, 0]$.
Next we compute the highest weight vector for $[0,0,0]$ with highest weight $[3, \frac{12}{7}, \frac{9}{7}]$. Set
\begin{alignat}{1}
& b^{1}=\omega^{1}_{-2}v^{3},\,\,\,\,\, b^{2}=\omega^{1}_{-1}\omega^{1}_{0}v^{3},\,\,\,\,\,
 b^{3}=\omega^{1}_{-1}\omega^{2}_{0}v^{3},\,\,\,\,\, b^{4}=\omega^{1}_{0}\omega^{1}_{0}\omega^{2}_{0}v^{3},
 \,\,\,\,\, b^{5}=\omega^{1}_{-1}\omega^{3}_{0}v^{3}, \nonumber \\
& b^{6}=\omega^{1}_{0}\omega^{1}_{0}\omega^{3}_{0}v^{3},\,\,\,\,\, b^{7}=\omega^{1}_{0}\omega^{2}_{-1}v^{3},\,\,\,\,\,
 b^{8}=\omega^{1}_{0}\omega^{2}_{0}\omega^{2}_{0}v^{3},\,\,\,\,\, b^{9}=\omega^{1}_{0}\omega^{2}_{0}\omega^{3}_{0}v^{3},
 \,\,\,\,\, b^{10}=\omega^{1}_{0}\omega^{3}_{-1}v^{3}, \nonumber \\
& b^{11}=\omega^{1}_{0}\omega^{3}_{0}\omega^{3}_{0}v^{3},\,\,\,\,\, b^{12}=\omega^{2}_{-2}v^{3},\,\,\,\,\,
 b^{13}=\omega^{2}_{-1}\omega^{2}_{0}v^{3},\,\,\,\,\, b^{14}=\omega^{2}_{0}\omega^{2}_{0}\omega^{2}_{0}v^{3},
 \,\,\,\,\, b^{15}=\omega^{2}_{-1}\omega^{3}_{0}v^{3}, \nonumber \\
& b^{16}=\omega^{2}_{0}\omega^{2}_{0}\omega^{3}_{0}v^{3},\,\,\,\,\, b^{17}=\omega^{2}_{0}\omega^{3}_{-1}v^{3},\,\,\,\,\,
 b^{18}=\omega^{2}_{0}\omega^{3}_{0}\omega^{3}_{0}v^{3},\,\,\,\,\, b^{19}=\omega^{3}_{-2}v^{3},
 \,\,\,\,\, b^{20}=\omega^{3}_{-1}\omega^{3}_{0}v^{3}, \nonumber \\
& b^{21}=\omega^{1}_{-1}u^{2},\,\,\,\,\, b^{22}=\omega^{2}_{-1}u^{2},\,\,\,\,\,
 b^{23}=\omega^{2}_{0}\omega^{2}_{0}u^{2},\,\,\,\,\, b^{24}=\omega^{2}_{0}\omega^{3}_{0}u^{2},
 \,\,\,\,\, b^{25}=\omega^{3}_{-1}u^{2}, \nonumber \\
& b^{26}=\omega^{3}_{0}\omega^{3}_{0}u^{2},\,\,\,\,\, b^{27}=v^{0}_{-2}v^{3}.\nonumber
\end{alignat}
By a similar method used in computations of $\langle t^{i}, \, t^{j}\rangle$, $1\leq i, j\leq 3$, we can calculate all  $\langle b^{i}, \, b^{j}\rangle$, $1\leq i, j\leq 27$. It is also
 a routine work to check that $\mathrm{det}\left(\langle b^{i},\, b^{j}\rangle\right)_{1\leq i, j\leq 27}\not=0$. Therefore $\{b^{i}, i=1,\ldots, 27\} $ is a linearly independent set. One can also see that
\begin{alignat}{1}
u^{4}=& -\frac{55}{486}\omega^{1}_{-2}v^{3}+\frac{5}{81}\omega^{1}_{-1}\omega^{1}_{0}v^{3}
-\frac{35}{18}\omega^{1}_{-1}\omega^{2}_{0}v^{3}+\frac{35}{36}\omega^{1}_{0}\omega^{1}_{0}\omega^{2}_{0}v^{3}\nonumber\\
& +\frac{7}{18}\omega^{1}_{0}\omega^{2}_{-1}v^{3}-\frac{49}{36}\omega^{1}_{0}\omega^{2}_{0}\omega^{2}_{0}v^{3}
+\frac{17}{27}\omega^{2}_{-2}v^{3}+\frac{7}{9}\omega^{2}_{-1}\omega^{2}_{0}v^{3}-
\frac{49}{27}\omega^{2}_{0}\omega^{2}_{0}\omega^{2}_{0}v^{3}\nonumber\\
& -\frac{5}{3}\omega^{1}_{-1}u^{2}-\frac{7}{99}\omega^{2}_{-1}u^{2}
-\frac{49}{132}\omega^{2}_{0}\omega^{2}_{0}u^{2}+v^{0}_{-2}v^{3}\nonumber
\end{alignat}
is a non-zero highest weight vector for $[0, 0, 0]$ with highest weight $[3, \frac{12}{7}, \frac{9}{7}]$. Next we show that $m$ is nonzero.
Set
\begin{alignat}{1}
& s^{i}=a^{i}, i=1,\ldots, 15,\,\,\,\,\, s^{16}=u^{3},\nonumber\\
& s^{j}=b^{j-16}, j=17,\ldots, 42,\,\,\,\,\, s^{43}=u^{4},\nonumber\\
& s^{44}=v^{2}_{-2}v^{3}.\nonumber
\end{alignat}
Then $\left(\mathcal{V}\otimes L\left(\frac{25}{28}, 0\right)\oplus \mathcal{V}\left(\frac{5}{7}\right)\otimes L\left(\frac{25}{28}, \frac{9}{7}\right)\right)^{(\tau_{e}, -)}_{5}$ is spanned by $s^{i}, i=1,\ldots, 43$.
By a similar method used in computations of $\langle t^{i}, \, t^{j}\rangle$, $1\leq i, j\leq 3$, we can calculate all  $\langle s^{i}, \, s^{j}\rangle$, $1\leq i, j\leq 44$. It is also
 a routine work to check that $\mathrm{det}\left(\langle s^{i},\, s^{j}\rangle\right)_{1\leq i, j\leq 44}\not=0$. Therefore $m$ must be nonzero. One
can also see that
\begin{alignat}{1}
u^{5}=\gamma^{1}+\gamma^{2}+v^{2}_{-2}v^{3}\nonumber
\end{alignat}
where
\begin{alignat}{1}
\gamma^{1}=& \frac{275}{486}\omega^{1}_{-2}v^{1}-\frac{25}{81}\omega^{1}_{-1}\omega^{1}_{0}v^{1}
+\frac{25}{72}\omega^{1}_{-1}\omega^{2}_{0}v^{1}-\frac{25}{144}\omega^{1}_{0}\omega^{1}_{0}\omega^{2}_{0}v^{1}
+\frac{5}{414}\omega^{1}_{0}\omega^{2}_{-1}v^{1}+\frac{65}{1656}\omega^{1}_{0}\omega^{2}_{0}\omega^{2}_{0}\nonumber\\
& +\frac{4}{5}\omega^{1}_{0}\omega^{3}_{-1}v^{1}-\frac{425}{11178}\omega^{2}_{-2}v^{1}
+\frac{325}{3726}\omega^{2}_{-1}\omega^{2}_{0}v^{1}+\frac{4}{5}\omega^{2}_{0}\omega^{3}_{-1}v^{1}
+\frac{4}{5}\omega^{3}_{-2}v^{1}+\frac{25}{39}\omega^{1}_{-1}u^{1}\nonumber\\
& +\frac{15}{26}\omega^{1}_{0}\omega^{1}_{0}u^{1}+\frac{5}{9}\omega^{2}_{-1}u^{1}
+\frac{72}{5}\omega^{3}_{-1}u^{1}+\frac{5}{782}u^{3}\nonumber
\end{alignat}
and
\begin{alignat}{1}
\gamma^{2}=& -\frac{605}{243}\omega^{1}_{-2}v^{3}+\frac{110}{81}\omega^{1}_{-1}\omega^{1}_{0}v^{3}
-\frac{55}{36}\omega^{1}_{-1}\omega^{2}_{0}v^{3}+\frac{55}{72}\omega^{1}_{0}\omega^{1}_{0}\omega^{2}_{0}v^{3}
-\frac{55}{36}\omega^{1}_{-1}\omega^{3}_{0}v^{3}+\frac{55}{72}\omega^{1}_{0}\omega^{1}_{0}\omega^{3}_{0}v^{3}\nonumber\\
& -\frac{13}{90}\omega^{1}_{0}\omega^{2}_{-1}v^{3}-\frac{17}{180}\omega^{1}_{0}\omega^{2}_{0}\omega^{2}_{0}v^{3}
-\frac{11}{36}\omega^{1}_{0}\omega^{2}_{0}\omega^{3}_{0}v^{3}-\frac{16}{23}\omega^{1}_{0}\omega^{3}_{-1}v^{3}+
\frac{20}{207}\omega^{1}_{0}\omega^{3}_{0}\omega^{3}_{0}v^{3}-\frac{13}{90}\omega^{2}_{-2}v^{3}\nonumber\\
& -\frac{13}{90}\omega^{2}_{-1}\omega^{2}_{0}v^{3}+\frac{1}{135}\omega^{2}_{0}\omega^{2}_{0}\omega^{2}_{0}v^{3}
-\frac{13}{90}\omega^{2}_{-1}\omega^{3}_{0}v^{3}-\frac{17}{180}\omega^{2}_{0}\omega^{2}_{0}\omega^{3}_{0}v^{3}
-\frac{16}{23}\omega^{2}_{0}\omega^{3}_{-1}v^{3}-\frac{188}{147}\omega^{3}_{-2}v^{3}\nonumber\\
& +\frac{20}{207}\omega^{2}_{0}\omega^{3}_{0}\omega^{3}_{0}v^{3}+\frac{40}{189}\omega^{3}_{-1}\omega^{3}_{0}v^{3}
+\frac{5}{6}\omega^{1}_{-1}u^{2}+\frac{11}{450}\omega^{2}_{-1}u^{2}
+\frac{41}{600}\omega^{2}_{0}\omega^{2}_{0}u^{2}+\frac{1}{8}\omega^{2}_{0}\omega^{3}_{0}u^{2}\nonumber\\
&+\frac{72}{253}\omega^{3}_{-1}u^{2}-\frac{10}{253}\omega^{3}_{0}\omega^{3}_{0}u^{2}-\frac{11}{5}u^{4}\nonumber
\end{alignat}
is a non-zero highest weight vector for $[0, 0, 0]$ with highest weight $[0, \frac{1}{7}, \frac{34}{7}]$. Thus the number $m$  must be $1$ following from  Theorem \ref{mirror}.
\end{proof}

\begin{corollary}
$\C \mathrm{VA}(e, f)$ is isomorphic to the $6A$-algebra constructed in \cite{LYY} as a VOA.
\end{corollary}
\begin{proof}
It is proved in \cite{DJY} a VOA with decomposition
\begin{alignat}{1}
V\cong \mathcal{V}\otimes L\left(\frac{25}{28}, 0\right)\oplus \mathcal{V}\left(\frac{1}{7}\right)\otimes L\left(\frac{25}{28}, \frac{34}{7}\right)\oplus \mathcal{V}\left(\frac{5}{7}\right)\otimes L\left(\frac{25}{28}, \frac{9}{7}\right).\nonumber
\end{alignat}
is unique and isomorphism to the $6A$-algebra in \cite{LYY}. So is $\C \mathrm{VA}(e, f)$.
\end{proof}
\subsection{Other cases}
Since the calculations are similar with what we did in the case $\langle e, f \rangle=\frac{5}{2^{10}}$, we only list some results and omit the details of the proof in this subsection. The following lemma is needed later \cite{DZh}:
\begin{lemma}\label{lemmaother}
For each $\mathcal{U}_{nX}, nX=1A, 2A, 2B, 4B$ constructed in \cite{LYY}, the VOA structure on $\mathcal{U}_{nX}$ is unique.
\end{lemma}
\subsubsection{The case $\langle e, f \rangle=\frac{1}{4}$}
In this case $e=f$. So $\C\mathrm{VA}(e, f)$ is generated by $e$. Thus $\C\mathrm{VA}(e, f)\cong L\left(\frac{1}{2},0\right)\cong \mathcal{U}_{1A}$ as VOAs.

\subsubsection{The case $\langle e, f \rangle=\frac{1}{2^5}$}
In this case $\G$ has a basis $e, f, \alpha(e, f)$. By \cite{S}, $\tau_e=\tau_f=\mathrm{id}$ on $\C\mathrm{VA}(e, f)$. So $e$ and $f$ are Ising vectors of $\sigma$ type. By \cite{JLY} and \cite{LYY}
, we have $\C\mathrm{VA}(e, f)\cong \mathcal{U}_{2A}$ as VOAs.

\subsubsection{The case $\langle e, f \rangle=0$}
In this case $\G$ has a basis $\{e, f\}$. By \cite{S}, $\tau_e=\tau_f=\mathrm{id}$ on $\C\mathrm{VA}(e, f)$. So $e$ and $f$ are Ising vectors of $\sigma$ type. By \cite{JLY} and \cite{LYY}
, we have $\C\mathrm{VA}(e, f)\cong \mathcal{U}_{2B}$ as VOAs.

\subsubsection{The case $\langle e, f \rangle=\frac{1}{2^8}$ with $(\tau_{e}\tau_{f})^4=1$}
In this case $\C\G$ has a basis $\{e, e^{\tau_{f}}, f, f^{\tau_{e}},  \alpha(e, f)\}$. Let $\omega$ be the Virasoro vector of $\C\mathrm{VA}(e, f)$. Let
\begin{alignat}{1}
& \omega^{1}=e, \nonumber\\
& \omega^{2}=\frac{1}{5}\left(3e+8e^{\tau_{f}}+4f+4f^{\tau_{e}}+2^{7}\alpha(e, f)\right),\nonumber\\
& \omega^{3}=-\frac{1}{5}\left(e+e^{\tau_{f}}-3f-3f^{\tau_{e}}+32\alpha(e, f)\right).\nonumber
\end{alignat}
Then $\omega^{1}$, $\omega^{2}$, $\omega^{3}$ are mutually orthogonal conformal vectors whose central charges are $\frac{1}{2}$, $\frac{7}{10}$, $\frac{7}{10}$ respectively and $\omega=\omega^{1}+\omega^{2}+\omega^{3}$ is  the Virasoro vector of $\C \mathrm{VA}(e, f).$
 We denote irreducible $L\left(\frac{1}{2},0\right)\otimes L\left(\frac{7}{10},0\right)\otimes L\left(\frac{7}{10},0\right)$-module $L\left(\frac{1}{2},h_{1}\right)\otimes L\left(\frac{7}{10},h_{2}\right)\otimes L\left(\frac{7}{10},h_{3}\right)$ as $[h_{1}, h_{2}, h_{3}]$. Then the following $L\left(\frac{1}{2},0\right)\otimes L\left(\frac{7}{10},0\right)\otimes L\left(\frac{7}{10},0\right)$-modules have integral weights:
\begin{alignat}{1}
& [0, 0, 0],\,\,\,\, [\frac{1}{2}, \frac{3}{2}, 0], \,\,\,\, [\frac{1}{16}, \frac{7}{16}, \frac{3}{2}], \nonumber\\
& [0, \frac{3}{2}, \frac{3}{2}],\,\,\,\, [\frac{1}{2}, 0, \frac{3}{2}], \,\,\,\, [\frac{1}{16}, \frac{3}{2}, \frac{7}{16}].\nonumber
\end{alignat}
 One can easily get that $[0, 0, 0], [\frac{1}{2}, \frac{3}{2}, 0],
 [0, \frac{3}{2}, \frac{3}{2}],\,\,\,\, [\frac{1}{2}, 0, \frac{3}{2}]$ are simple current modules for the VOA $L\left(\frac{1}{2},0\right)\otimes L\left(\frac{7}{10},0\right)\otimes L\left(\frac{7}{10},0\right).$

Let
\begin{alignat}{1}
& v^{1}=e+f+f^{\tau_{e}}+32\alpha(e, f), \nonumber\\
& v^{2}=f-f^{\tau_{e}}.\nonumber
\end{alignat}
Then we have
\begin{alignat}{1}
& \omega^{1}v^{1}=\frac{1}{2}v^{1},\,\,\,\, \omega^{2}v^{1}=\frac{3}{2}v^{1},\,\,\,\, \omega^{3}v^{1}=0, \nonumber\\
& \omega^{1}v^{2}=\frac{1}{16}v^{2},\,\,\,\, \omega^{2}v^{2}=\frac{7}{16}v^{2},\,\,\,\, \omega^{3}v^{2}=\frac{3}{2}v^{2}.\label{4Bv1v2}
\end{alignat}
Since the multiplicity of a simple current module in an extension VOA is $0$ or $1$ \cite{DMZ}, we have
\begin{alignat}{1}
\C\mathrm{VA}(e, f)\cong [0, 0, 0]\oplus [\frac{1}{2}, \frac{3}{2}, 0]\oplus a\cdot [0, \frac{3}{2}, \frac{3}{2}] \oplus b\cdot [\frac{1}{2}, 0, \frac{3}{2}]
\oplus m\cdot [\frac{1}{16}, \frac{7}{16}, \frac{3}{2}]\oplus n\cdot [\frac{1}{16}, \frac{3}{2}, \frac{7}{16}]\nonumber
\end{alignat}
as $L\left(\frac{1}{2},0\right)\otimes L\left(\frac{7}{10},0\right)\otimes L\left(\frac{7}{10},0\right)$-modules, where $a$ and $b$ are $0$ or $1$, $m,n \in \Z_{+}$ and $m\neq 0.$

Using the $\tau$-involution  and $\sigma$-involution defined in \S \ref{subsec:The-unitary-series}, 
 we apply Theorem \ref{quantum dimension and orbifold module} to $G=\langle\tau_{\omega^1}\rangle$, $\langle\tau_{\omega^2}\rangle$, $\langle\sigma_{\omega^3}\rangle$ and $V= \C VA(e,f).$ Then we immediately get
\begin{equation*}
\left\{
\begin{array}{cccccccc}
2+a+b & =  & 2m+2n, \\[0.3cm]
2+a+b+2n & =  & 2m,\\[0.3cm]
 2 & =  & a+b+2m.\\[0.3cm]
\end{array}
\right.
\end{equation*}

So $n=a=b=0$ and $m=1,$ i.e.
\begin{alignat}{1}
\C\mathrm{VA}(e, f)\cong [0, 0, 0]\oplus [\frac{1}{2}, \frac{3}{2}, 0] \oplus [\frac{1}{16}, \frac{7}{16}, \frac{3}{2}].\nonumber
\end{alignat}
By \cite{LYY} and Lemma \ref{lemmaother}, a VOA with such a decomposition is unique. Thus we have $\C\mathrm{VA}(e, f)\cong \mathcal{U}_{4B}$ as VOAs.


\section*{Acknowledgments}
The author Wen Zheng thanks Chongying, Dong for his discussions and helpful comments.

\end{document}